\documentclass[11pt,a4paper]{amsart}
\usepackage[T1]{fontenc}
\usepackage{lmodern}
\usepackage{amsmath,amssymb,amsthm}
\usepackage{float}
\usepackage[hidelinks]{hyperref}

\newtheorem{theorem}{Theorem}[section]
\newtheorem{proposition}[theorem]{Proposition}
\newtheorem{lemma}[theorem]{Lemma}

\newtheorem*{conjecture}{Conjecture}
\numberwithin{equation}{section}
\newcommand{\N}{\mathbb N}
\newcommand{\Pp}{\mathbb P}
\newcommand{\Pc}{\mathcal P}

\newcommand{\Dc}{\mathcal D}
\newcommand{\Ec}{\mathcal E}
\newcommand{\Gc}{\mathcal G}
\newcommand{\meas}{\operatorname{meas}}

\title[A counterexample to Beurling's conjecture]{A Counterexample to Beurling's Conjecture}
\author{Volodymyr Pavlenkov}
\address{Department of Mathematics, University of York,
Heslington, York, YO10 5DD, England}
\email{volodymr.pavlenkov@york.ac.uk}
\subjclass[2020]{Primary 11N80; Secondary 11N05, 11M26}
\keywords{Beurling generalized primes, Beurling's conjecture, Cram\'er's conjecture, primes in short intervals, exceptional sets}

\begin{document}
\begin{abstract}
Beurling conjectured that a system of generalized primes whose associated integer-counting function satisfies $N(x)=x+o(\log x)$ must coincide with the system of ordinary primes. We construct a nonclassical generalized prime system for which $N(x)=x+O((\log x)^\beta)$, where $\beta=\log 2/\log(201/100)<1$. The system is obtained by deleting the primes $2,3,5$ and replacing an infinite sequence of ordinary primes by their square roots. The construction of this sequence uses the exceptional-set estimates for the prime number theorem in short intervals established by Gafni and Tao. In particular, Beurling's conjecture is false. We propose a revised conjecture with the error term $o((\log x)^{\log 2/\log 3})$. Under Cram\'er's conjecture, we construct nonclassical systems with error $O((\log x)^{\log 2/\log 3})$, and show that the implied constant can be made arbitrarily small.
\end{abstract}
\maketitle

\section{Introduction}

Let $\Pc=\{q_1,q_2,\ldots\}$ be a sequence of real numbers satisfying
\[
  1<q_1<q_2<\cdots,\qquad q_j\to\infty.
\]
The elements of $\Pc$ are called Beurling generalized primes. The associated generalized integers are the formal products
\[
  \prod_{j\geq1}q_j^{\nu_j},\qquad
  \nu_j\in\mathbb Z_{\geq0},\quad \nu_j=0\text{ for all but finitely many }j.
\]
We denote the corresponding integer-counting function by $N_{\Pc}(x)$, counting formal products of numerical value at most $x$, with multiplicities. We also put
\[
  \pi_{\Pc}(x)=\#\{q\in\Pc:q\leq x\}.
\]
Let $\Pp$ denote the set of ordinary primes and $\pi(x)$ the ordinary prime-counting function. By the fundamental theorem of arithmetic, $N_{\Pp}(x)=\lfloor x\rfloor$.

The following conjecture of Beurling is recorded in Erd\H{o}s' survey on combinatorial number theory \cite[p.~103, (12)]{Erdos}; see also \cite[Conjecture~1.2]{Olofsson} and the commentary to \cite[Problem~951]{Bloom}.

\begin{conjecture}[Beurling]
If
\begin{equation}\label{eq:conj}
  N_{\Pc}(x)=x+o(\log x),\qquad x\to\infty,
\end{equation}
then $\Pc=\Pp$.
\end{conjecture}

Erd\H{o}s formulates the conjecture in terms of the set of numerical values of the products, with the strict inequality in the counting function. The counting convention used here is that of modern Beurling theory. In the system constructed in this paper distinct formal products have distinct numerical values, so the counterexample applies also to the distinct-value formulation with unique factorization. The two conventions concerning the endpoint $x$ then differ by at most one. We do not impose the unit separation condition occurring in the related question of Shapiro \cite[p.~103, (13)]{Erdos}.

Olofsson \cite[Theorem~1.3]{Olofsson} proved that for every $c>0$ there is a generalized prime system $\Pc\ne\Pp$ satisfying $|N_{\Pc}(x)-\lfloor x\rfloor|<c\log x$ for all $x>1$. Thus, if \eqref{eq:conj} implied $\Pc=\Pp$, the order of the error term in that conjecture would be optimal. We show that the logarithmic order can in fact be improved for a nonclassical system.

\begin{theorem}\label{thm:main}
There exists a sequence of ordinary primes $5<r_1<r_2<\cdots$ such that the generalized prime system
\begin{equation}\label{eq:system}
  \Pc^*=\big(\Pp\setminus\{2,3,5,r_1,r_2,\ldots\}\big)
  \cup\{\sqrt{r_1},\sqrt{r_2},\ldots\}
\end{equation}
satisfies
\begin{equation}\label{eq:main}
  N_{\Pc^*}(x)=x+O\big((\log x)^\beta\big),\qquad x\to\infty,
\end{equation}
where
\[
  \beta=\frac{\log 2}{\log(201/100)}=0.9928559\ldots<1.
\]
Moreover, distinct formal products of elements of $\Pc^*$ have distinct numerical values and
\begin{equation}\label{eq:diff}
  \liminf_{x\to\infty}(\pi_{\Pc^*}(x)-\pi(x))=-3,
  \qquad
  \limsup_{x\to\infty}(\pi_{\Pc^*}(x)-\pi(x))=-2.
\end{equation}
\end{theorem}

Since $\Pc^*\ne\Pp$ and $\beta<1$, Theorem~\ref{thm:main} disproves Beurling's conjecture. The additional relations \eqref{eq:diff} show that a counterexample can be constructed even under the assumption that the generalized prime-counting function differs from the ordinary prime-counting function by a bounded quantity and is strictly smaller for all sufficiently large $x$.

The construction reduces to finding a sufficiently sparse sequence of ordinary primes satisfying the identity
\begin{equation}\label{eq:productintro}
  \prod_{j=1}^{\infty}(1+r_j^{-1/2})=\frac{15}{4}.
\end{equation}
The square-root construction was introduced in the companion note \cite{Pavlenkov}, where it is shown that Cram\'er's conjecture implies the failure of Beurling's conjecture and, more generally, that a uniform bound $p_{n+1}-p_n=O(p_n^\theta)$ with $\theta<1/2$ suffices. In the present paper we replace that hypothesis by an unconditional estimate on the measure of exceptional intervals for the prime number theorem in short intervals. The main step is to select a prime at each stage of the construction in such a way that the next stage admits a sufficiently large number of admissible primes.

Theorem~\ref{thm:main} also raises the question of the order of the error term for which rigidity can still be expected. Put $\beta_* = \log 2/\log 3$. In Section~5 we refine the conditional estimate of \cite{Pavlenkov} to the endpoint $\beta_*$ and propose a revised conjecture: the relation $N_{\Pc}(x)=x+o((\log x)^{\beta_*})$ implies $\Pc=\Pp$. Under Cram\'er's conjecture no larger exponent can satisfy such a rigidity statement. We also obtain a conditional analogue of Olofsson's result at this scale.

In Section~2 we derive the integer-counting estimate from the normalization \eqref{eq:productintro} and a growth condition on $(r_j)$. In Section~3 we establish an auxiliary exceptional-set estimate from the results of Gafni and Tao \cite{GT}. In Section~4 we construct the required sequence and complete the proof of Theorem~\ref{thm:main}. Section~5 is devoted to the revised conjecture and its conditional sharpness.

\section{The square-root construction}

\begin{proposition}\label{prop:count}
Let $(r_j)$ be a strictly increasing sequence of ordinary primes greater than $5$ satisfying \eqref{eq:productintro}. Suppose that, for some $\lambda>2$ and $c>0$,
\begin{equation}\label{eq:growth}
  r_{j+1}\geq c r_j^\lambda
\end{equation}
for all sufficiently large $j$. Then the system \eqref{eq:system} has distinct numerical values for distinct formal products, satisfies \eqref{eq:diff}, and
\begin{equation}\label{eq:countconcl}
  N_{\Pc^*}(x)=x+O\big((\log x)^{\log 2/\log\lambda}\big).
\end{equation}
\end{proposition}

\begin{proof}
Each ordinary prime greater than $5$ is either an element of $\Pc^*$ or the square of an element of $\Pc^*$. Squaring a formal product of elements of $\Pc^*$ gives an ordinary integer: the exponent of an unmodified prime is twice its formal exponent, and the exponent of $r_j$ is the formal exponent of $\sqrt{r_j}$. By uniqueness of ordinary prime factorization, two formal products with equal numerical values must have identical formal exponents.

Put
\[
\begin{gathered}
  \Dc=\left\{\prod_{j\in J}r_j:J\subset\N,\ |J|<\infty\right\},\\
  B(y)=\#\{m\in\N:m\leq y,\ (m,30)=1\},\quad y\geq0,
\end{gathered}
\]
where the empty product is $1$. Decomposing each exponent of $\sqrt{r_j}$ into an even part and a remainder in $\{0,1\}$, every generalized integer has a unique representation $m\sqrt d$ with $(m,30)=1$ and $d\in\Dc$. Consequently,
\begin{equation}\label{eq:exactcount}
  N_{\Pc^*}(x)=\sum_{\substack{d\in\Dc\\d\leq x^2}}B(x/\sqrt d).
\end{equation}
By inclusion--exclusion,
\[
  B(y)=\sum_{e\mid30}\mu(e)\left\lfloor\frac{y}{e}\right\rfloor
      =\frac4{15}y+O(1),\qquad y\geq0,
\]
where the implied constant is uniform in $y$.

For $\sqrt{r_k}\leq x<\sqrt{r_{k+1}}$, only subsets of $\{1,\ldots,k\}$ can contribute to \eqref{eq:exactcount}. Since $B(y)=0$ for $0\leq y<1$, we can include all $2^k$ such subsets, including those for which $d>x^2$. Hence
\begin{equation}\label{eq:truncated}
  N_{\Pc^*}(x)=\frac4{15}x\prod_{j=1}^{k}(1+r_j^{-1/2})+O(2^k).
\end{equation}
By \eqref{eq:growth}, $r_{j+1}\geq4r_j$ eventually. Thus
\[
  S_k:=\sum_{j>k}r_j^{-1/2}\leq2r_{k+1}^{-1/2}
\]
for all sufficiently large $k$. From \eqref{eq:productintro},
\[
  0\leq 1-\frac4{15}\prod_{j=1}^{k}(1+r_j^{-1/2})
  =1-\prod_{j>k}(1+r_j^{-1/2})^{-1}\leq S_k.
\]
Here the last inequality follows from $1-e^{-v}\leq v$ for $v\geq0$ and $\log(1+w)\leq w$ for $w\geq0$. Since $x<\sqrt{r_{k+1}}$, \eqref{eq:truncated} now gives
\begin{equation}\label{eq:2k}
  N_{\Pc^*}(x)=x+O(2^k).
\end{equation}

Set $u_j=\log r_j$. By \eqref{eq:growth}, $u_{j+1}\geq\lambda u_j+\log c$ eventually. Choose $j_0$ such that this inequality holds for $j\geq j_0$ and $u_{j_0}+\log c/(\lambda-1)>0$. Iteration gives $u_j\gg\lambda^j$. Consequently,
\[
  2^k\ll (\log r_k)^{\log 2/\log\lambda}
       \ll (\log x)^{\log 2/\log\lambda}.
\]
Together with \eqref{eq:2k}, this proves \eqref{eq:countconcl}.

Finally, putting $A(x)=\#\{j:r_j\leq x\}$, we obtain from \eqref{eq:system}, for $x\geq5$,
\begin{equation}\label{eq:delta}
  \pi_{\Pc^*}(x)-\pi(x)=-3+A(x^2)-A(x).
\end{equation}
The difference $A(x^2)-A(x)$ counts the intervals $[\sqrt{r_j},r_j)$ containing $x$. Since $\lambda>2$, \eqref{eq:growth} implies $\sqrt{r_{j+1}}>r_j$ for all sufficiently large $j$. The intervals are therefore eventually pairwise disjoint and have nonempty gaps between consecutive intervals. For all sufficiently large $x$, the right-hand side of \eqref{eq:delta} belongs to $\{-3,-2\}$, and both values occur arbitrarily far out. This establishes \eqref{eq:diff}.
\end{proof}

\section{An estimate for exceptional intervals}

Throughout Sections~3 and~4 put
\begin{equation}\label{eq:parameters}
  a=\frac{99}{200},\qquad \lambda=3-2a=\frac{201}{100}.
\end{equation}
Let $\psi(x)=\sum_{n\leq x}\Lambda(n)$ be the classical Chebyshev function, and define
\[
  \Ec=\left\{y\geq1:\left|\psi(y+y^a)-\psi(y)-y^a\right|\geq\frac12 y^a\right\}.
\]

\begin{lemma}\label{lem:exception}
For $X\to\infty$,
\begin{equation}\label{eq:exception}
  \meas(\Ec\cap[X,2X])\ll X^{37/50}.
\end{equation}
\end{lemma}

\begin{proof}
We use the notation $A(\sigma)$ and $A^*(\sigma)$ of Gafni and Tao \cite[Section~1]{GT} for the zero-density and zero-density additive-energy exponents, respectively. In particular, if $N(\sigma,T)$ denotes the number of zeros $\rho$ of the Riemann zeta function with $\Re\rho\geq\sigma$ and $|\Im\rho|\leq T$, counted with multiplicities, then $A(\sigma)$ is the infimum of the exponents $A$ for which $N(\sigma,T)\ll_{\sigma,\varepsilon}T^{(1-\sigma)A+\varepsilon}$ for each $\varepsilon>0$. The exponent $A^*(\sigma)$ is defined analogously using the number of quadruples of these zeros with
\[
  |\Im(\rho_1+\rho_2-\rho_3-\rho_4)|\leq1.
\]
In particular, $H(\sigma):=(1-\sigma)A^*(\sigma)$ is nonincreasing, by inclusion of the sets of quadruples as $\sigma$ increases. Let $\mu(\theta)$ denote the exceptional-set exponent in \cite[Definition~1.1]{GT}.

By \cite[Theorem~1.3]{GT},
\begin{equation}\label{eq:GT}
\begin{split}
  \mu(a)\leq\max\bigg\{1-a,\ \inf_{\varepsilon>0}\!
  \sup_{\substack{1/2<\sigma<1\\ A(\sigma)\geq(1-a)^{-1}-\varepsilon}}
  \min\{M_2(\sigma),M_4(\sigma)\}\bigg\},
\end{split}
\end{equation}
where
\[
  M_2(\sigma)=(1-a)(1-\sigma)A(\sigma)+2\sigma-1,
  \qquad M_4(\sigma)=(1-a)H(\sigma)+4\sigma-3.
\]
For $4/5\leq\sigma<7/8$, Table~1 of \cite{GT} gives $A(\sigma)\leq3/(2\sigma)\leq15/8$. For $7/8\leq\sigma<1$, the Ingham--Huxley estimate $A(\sigma)\leq3/(3\sigma-1)$, recalled in \cite[Section~3]{GT}, gives $A(\sigma)\leq24/13<15/8$. Since $(1-a)^{-1}=200/101>15/8$, for sufficiently small $\varepsilon>0$ the supremum in \eqref{eq:GT} is restricted to $1/2<\sigma<4/5$.

We next give a rational upper bound for $M_4$ on this entire interval. For $i=100,\ldots,159$ put $s_i=i/200$. Let $H_i$ be the upper bound for $H(s_i)$ obtained from Table~2 of \cite{GT}. The expressions needed at these rational grid points are listed in Table~\ref{tab:density}; they are obtained by multiplying the corresponding bound for $A^*(s_i)$ by $1-s_i$. By monotonicity of $H$, for $s_i\leq\sigma\leq s_{i+1}$,
\begin{equation}\label{eq:grid}
  M_4(\sigma)\leq \frac{101}{200}H_i+4s_{i+1}-3=:V_i.
\end{equation}
\begin{table}[H]
\centering
\begin{tabular}{c|c}
$i$ & $H_i$ (with $s=s_i$) \\ \hline
$100,\ldots,133$ & $\dfrac{10-11s}{2-s}$ \\[5pt]
$134,\ldots,140$ & $\dfrac{18-19s}{4-2s}$ \\[5pt]
$141,\ldots,145$; $148,\ldots,152$ & $\dfrac{5(18-19s)}{2(5s+3)}$ \\[5pt]
$146$ & $\dfrac{2(45-44s)}{2s+15}$ \\[5pt]
$147$ & $\dfrac{457-546s}{2(61-58s)}$ \\[5pt]
$153$ & $\dfrac{3(18-19s)}{4(4s-1)}$ \\[5pt]
$154$ & $\dfrac{5(18-19s)}{2(13s-3)}$ \\[5pt]
$155,\ldots,159$ & $\dfrac{4(10-9s)}{5(4s-1)}$
\end{tabular}
\caption{Upper bounds at the grid points from \cite[Table~2]{GT}.}\label{tab:density}
\end{table}
Substituting $s_i=i/200$ into \eqref{eq:grid} and comparing the resulting rational numbers gives the maxima in Table~\ref{tab:verification}. In particular,
\[
  \max_{100\leq i\leq159}V_i=\frac{3811}{5200}<\frac{37}{50}.
\]
The finite calculation in Table~\ref{tab:verification} uses only the bounds at the listed rational grid points and monotonicity of $H$ between them; no interpolation of a numerical plot is required.
\begin{table}[H]
\centering
\begin{tabular}{c|c|c}
Range of $i$ & Maximizing $i$ & $\max V_i$ \\ \hline
$100,\ldots,109$ & $109$ & $11447/19400$ \\
$110,\ldots,119$ & $119$ & $36071/56200$ \\
$120,\ldots,129$ & $129$ & $37001/54200$ \\
$130,\ldots,139$ & $139$ & $75979/104400$ \\
$140,\ldots,149$ & $140$ & $3811/5200$ \\
$150,\ldots,159$ & $159$ & $79269/109000$
\end{tabular}
\caption{Exact rational maxima in \eqref{eq:grid}.}\label{tab:verification}
\end{table}

Since $1-a=101/200<3811/5200$, \eqref{eq:GT} implies
\[
  \mu(a)\leq\frac{3811}{5200}<\frac{37}{50}.
\]
Applying the definition of the exceptional-set exponent with deviation parameter $1/2$ proves \eqref{eq:exception}.
\end{proof}

\section{Construction of the sequence of ordinary primes}

Put
\begin{equation}\label{eq:f}
  f(t)=\log(1+t^{-1/2}),\qquad t>0.
\end{equation}
The function $f$ is strictly decreasing from $+\infty$ to $0$, and
\begin{equation}\label{eq:fproperties}
  f^{-1}(v)=(e^v-1)^{-2}\sim v^{-2}\quad(v\downarrow0),
  \qquad f'(t)=-\frac1{2t(\sqrt t+1)}.
\end{equation}
We shall construct a strictly increasing sequence of ordinary primes satisfying
\begin{equation}\label{eq:sum}
  \sum_{j=1}^\infty f(r_j)=\log(15/4),
  \qquad r_{j+1}\asymp r_j^\lambda
\end{equation}
for all sufficiently large $j$.

Define the set
\begin{equation}\label{eq:good}
  \Gc=\left\{t\geq2:\pi(t+3t^a)-\pi(t+t^a)\geq\frac{t^a}{10\log t}\right\}.
\end{equation}
For $t>0$ and $r>t$ put
\begin{equation}\label{eq:T}
  T_t(r)=f^{-1}(f(t)-f(r)).
\end{equation}
The quantity $T_t(r)$ will determine the location of the primes admissible at the next step of the construction.

\begin{lemma}\label{lem:bad}
If $u\notin\Gc$ is sufficiently large, then
\begin{equation}\label{eq:badinterval}
  [u+\tfrac{11}{10}u^a,u+\tfrac65u^a]\subset\Ec.
\end{equation}
\end{lemma}
\begin{proof}
For $y$ in the interval on the left of \eqref{eq:badinterval}, $y^a=(1+o(1))u^a$ uniformly as $u\to\infty$. In particular, $(y,y+y^a]\subset(u+u^a,u+3u^a]$ for large $u$. By \eqref{eq:good}, the contribution of ordinary primes to $\psi(y+y^a)-\psi(y)$ is at most
\[
  \frac{u^a}{10\log u}\log(u+3u^a)=(1/10+o(1))u^a.
\]
If two distinct $m$th powers $h^m$ and $k^m$, $2\leq m$, lie in this interval, then, taking $k>h$, the mean value theorem gives $k^m-h^m\geq m h^{m-1}\geq 2u^{1/2}$. Since $a<1/2$, there is at most one such power in $(u+u^a,u+3u^a]$ for each $m$, for large $u$. There are $O(\log u)$ possible exponents $m$, and each prime power contributes at most $O(\log u)$ to $\psi$. The contribution of the higher prime powers is therefore $O((\log u)^2)=o(u^a)$. Consequently,
\[
  \psi(y+y^a)-\psi(y)\leq(1/10+o(1))u^a<\tfrac12 y^a.
\]
By the definition of $\Ec$, this proves \eqref{eq:badinterval}.
\end{proof}

\begin{lemma}\label{lem:selection}
Let $b\in\{a,3/5\}$. Suppose that $t$ is sufficiently large and, for a fixed $c_0>0$,
\begin{equation}\label{eq:candidates}
  \pi(t+3t^b)-\pi(t+t^b)\geq c_0\frac{t^b}{\log t}.
\end{equation}
Then there exists an ordinary prime $r\in(t+t^b,t+3t^b]$ such that $T_t(r)\in\Gc$.
\end{lemma}

\begin{proof}
For $r\in[t+t^b,t+3t^b]$, the mean value theorem and \eqref{eq:fproperties} give, uniformly in $r$,
\[
  f(t)-f(r)\asymp t^{b-3/2},\qquad T_t(r)\asymp t^{3-2b}.
\]
Put $Y=t^{3-2b}$. There exist constants $c_1,c_2>0$ such that all the values $T_t(r)$ lie in $[c_1Y,c_2Y]$ for large $t$. Differentiating $f(T_t(r))=f(t)-f(r)$ and applying \eqref{eq:fproperties}, we obtain
\begin{equation}\label{eq:derivative}
  |T_t'(r)|=\frac{|f'(r)|}{|f'(T_t(r))|}\asymp t^{3-3b}.
\end{equation}
The map $r\mapsto T_t(r)$ is strictly decreasing. Since distinct ordinary primes differ by at least one, the mean value theorem and \eqref{eq:derivative} show that distinct candidate primes in \eqref{eq:candidates} give target values separated by $\gg t^{3-3b}$.

Suppose that none of the target values belongs to $\Gc$. By Lemma~\ref{lem:bad}, each target $u=T_t(r)$ then gives an interval of length $u^a/10\asymp Y^a$ contained in $\Ec$. Since
\[
  3-3b>a(3-2b)
\]
for both $b=a$ and $b=3/5$, the separation of the targets is $\gg t^{3-3b}\gg Y^a$ with ratio tending to infinity. The intervals in \eqref{eq:badinterval} associated with distinct candidate primes are thus pairwise disjoint for sufficiently large $t$. They all lie in $[c_3Y,c_4Y]$ for fixed positive constants $c_3,c_4$. By \eqref{eq:candidates},
\begin{equation}\label{eq:lower}
  \meas(\Ec\cap[c_3Y,c_4Y])\gg \frac{t^bY^a}{\log t}
  \asymp \frac{Y^{a+b/(3-2b)}}{\log Y}.
\end{equation}
On the other hand, covering $[c_3Y,c_4Y]$ by a bounded number of dyadic intervals and applying Lemma~\ref{lem:exception} gives
\begin{equation}\label{eq:upper}
  \meas(\Ec\cap[c_3Y,c_4Y])\ll Y^{37/50}.
\end{equation}
By \eqref{eq:parameters},
\[
  a+\frac{a}{3-2a}=\frac{9933}{13400}>\frac{37}{50},\qquad
  a+\frac{3/5}{3-6/5}=\frac{497}{600}>\frac{37}{50}.
\]
Thus \eqref{eq:lower} contradicts \eqref{eq:upper} for large $Y$. The assertion follows.
\end{proof}

\begin{proposition}\label{prop:sequence}
For every $C>0$ and $M>5$, there exists a strictly increasing sequence of ordinary primes $r_j>M$ such that
\[
  \sum_{j=1}^\infty f(r_j)=C,
  \qquad r_{j+1}\asymp r_j^{201/100}
\]
for all sufficiently large $j$.
\end{proposition}

\begin{proof}
We first record two consequences of Lemma~\ref{lem:selection}. If $t\in\Gc$ is sufficiently large, then \eqref{eq:candidates} holds for $b=a$ by definition, so an ordinary prime $r\in(t+t^a,t+3t^a]$ can be selected with $T_t(r)\in\Gc$.

For an arbitrary sufficiently large $t$, we apply the prime number theorem in short intervals of length $z^{3/5}$ (see \cite[Theorem~1.1 and the subsequent discussion]{GT}). Taking $z=t+t^{3/5}$, the interval $(z,z+z^{3/5}]$ is contained in $(t+t^{3/5},t+3t^{3/5}]$ for sufficiently large $t$. The higher prime powers contribute $O(t^{1/2}(\log t)^2)=o(t^{3/5})$ to the Chebyshev function in this interval. The prime number theorem in short intervals thus implies \eqref{eq:candidates} for $b=3/5$, with some absolute $c_0>0$. Lemma~\ref{lem:selection} again provides an ordinary prime $r$ for which $T_t(r)\in\Gc$. Moreover, in both cases $T_t(r)\asymp t^{3-2b}$ and $r\asymp t$, so $T_t(r)>r$ for all sufficiently large $t$.

Choose $R>M$ so large that these conclusions hold for every $t\geq R$ and $f(R)<C$. Enumerate the ordinary primes greater than $R$ as $s_1<s_2<\cdots$. Since $f(s)\asymp s^{-1/2}$ as $s\to\infty$, by the divergence of the sum of reciprocals of the ordinary primes we have $\sum_j f(s_j)=+\infty$. Let $m$ be the largest nonnegative integer satisfying
\[
  \sum_{j=1}^{m} f(s_j)<C.
\]
Our choice of $R$ gives $m\geq1$. Put $r_j=s_j$ for $1\leq j\leq m$, and set
\[
  \delta_m=C-\sum_{j=1}^{m}f(r_j),\qquad t_m=f^{-1}(\delta_m).
\]
By the definition of $m$, $0<\delta_m\leq f(s_{m+1})$, so $t_m\geq s_{m+1}>r_m$ and $t_m>R$.

By the second consequence of Lemma~\ref{lem:selection}, choose
\[
  r_{m+1}\in(t_m+t_m^{3/5},t_m+3t_m^{3/5}]
\]
to be an ordinary prime with $t_{m+1}:=T_{t_m}(r_{m+1})\in\Gc$. Then $t_{m+1}>r_{m+1}$. Suppose now that $r_k$ and $t_k\in\Gc$ have been constructed for some $k\geq m+1$, with $t_k>r_k$ and
\begin{equation}\label{eq:residual}
  f(t_k)=C-\sum_{j=1}^{k}f(r_j)>0.
\end{equation}
Using the first consequence of Lemma~\ref{lem:selection}, choose an ordinary prime
\begin{equation}\label{eq:rk}
  r_{k+1}\in(t_k+t_k^a,t_k+3t_k^a]
\end{equation}
such that $t_{k+1}:=T_{t_k}(r_{k+1})\in\Gc$. Then $t_{k+1}>r_{k+1}>t_k>r_k$ and, by \eqref{eq:T},
\[
  f(t_{k+1})=f(t_k)-f(r_{k+1})=C-\sum_{j=1}^{k+1}f(r_j)>0.
\]
This proves \eqref{eq:residual} inductively, with the initial case $k=m+1$ following from the definition of $t_{m+1}$. The selection may be made deterministic by always taking the least admissible prime.

By \eqref{eq:rk}, $r_{k+1}\asymp t_k$. The estimates preceding \eqref{eq:derivative}, with $b=a$, give $t_{k+1}\asymp t_k^\lambda$. Consequently,
\[
  r_{k+2}\asymp t_{k+1}\asymp t_k^\lambda\asymp r_{k+1}^\lambda,
  \qquad k\geq m+1.
\]
Since $(r_k)$ is a strictly increasing sequence of primes and $t_k>r_k$, we have $t_k\to\infty$. Passing to the limit in \eqref{eq:residual} and using $f(t_k)\to0$, we obtain $\sum_{j\geq1}f(r_j)=C$. This proves the proposition.
\end{proof}

\begin{proof}[Proof of Theorem~\ref{thm:main}]
Apply Proposition~\ref{prop:sequence} with $C=\log(15/4)$ and any $M>5$. By \eqref{eq:f}, the resulting sequence satisfies \eqref{eq:productintro}. Its growth estimate implies \eqref{eq:growth} with $\lambda=201/100>2$. Proposition~\ref{prop:count} establishes \eqref{eq:main}, \eqref{eq:diff} and uniqueness of numerical values of formal products. Since the prime $2$ does not belong to $\Pc^*$, we have $\Pc^*\ne\Pp$. Finally, $\log 2/\log(201/100)<1$, so \eqref{eq:main} implies \eqref{eq:conj}.
\end{proof}

\section{A revised conjecture}

In \cite{Pavlenkov} we proved that under Cram\'er's conjecture, for every $\beta>\log 2/\log 3$, there exists a nonclassical generalized prime system with integer-counting function $N_{\Pc}(x)=x+O((\log x)^\beta)$. We now establish the corresponding estimate at the endpoint and show that the constant in the error estimate can be made arbitrarily small. Put
\[
  \beta_* = \frac{\log 2}{\log 3}.
\]
Recall that Cram\'er's conjecture \cite{Cramer} asserts, in its upper-bound form, that the ordinary primes, enumerated in increasing order, satisfy
\begin{equation}\label{eq:cramer}
  p_{n+1}-p_n=O((\log p_n)^2).
\end{equation}

\begin{proposition}\label{prop:endpoint}
Assume Cram\'er's conjecture. Then for every $c>0$ there exists a generalized prime system $\Pc\ne\Pp$ with distinct numerical values for distinct formal products such that
\begin{equation}\label{eq:smallconstant}
  |N_{\Pc}(x)-x|\leq c(\log x)^{\beta_*}
\end{equation}
for all sufficiently large $x$.
\end{proposition}

\begin{proof}
Let $p$ be a sufficiently large ordinary prime. We shall construct a strictly increasing sequence $(r_j)$ of ordinary primes greater than $p$ such that
\begin{equation}\label{eq:pnormalization}
  \prod_{j=1}^\infty(1+r_j^{-1/2})=\frac{p}{p-1}.
\end{equation}
Let $f$ be the function in \eqref{eq:f} and let $P(t)$ denote the least ordinary prime strictly greater than $t$. By \eqref{eq:cramer},
\[
  0<P(t)-t\ll (\log t)^2,\qquad P(t)\sim t.
\]
By the mean value theorem and \eqref{eq:fproperties}, for all sufficiently large $t$,
\begin{equation}\label{eq:cramerstep}
  0<f(t)-f(P(t))\ll t^{-3/2}(\log t)^2.
\end{equation}
In particular, $f(t)-f(P(t))<f(t)/2$ for all sufficiently large $t$.

Set $\delta_0=\log(p/(p-1))$ and $t_0=f^{-1}(\delta_0)=(p-1)^2$. Define recursively, for $k\geq0$,
\[
  r_{k+1}=P(t_k),\qquad
  \delta_{k+1}=f(t_k)-f(r_{k+1}),\qquad
  t_{k+1}=f^{-1}(\delta_{k+1}).
\]
For $p$ sufficiently large, \eqref{eq:cramerstep} implies $0<\delta_{k+1}<\delta_k/2$ and $\delta_{k+1}<f(r_{k+1})$. Since $f$ is strictly decreasing, $t_{k+1}>r_{k+1}>t_k$. Thus the recursion is well defined, $r_1>p$, and $(r_j)$ is strictly increasing. Moreover,
\[
  \delta_k=\log\frac{p}{p-1}-\sum_{j=1}^k f(r_j)\longrightarrow0,
\]
which establishes \eqref{eq:pnormalization}.

By \eqref{eq:cramerstep}, the relation $r_{k+1}\sim t_k$, and the asymptotic relation for $f^{-1}$ in \eqref{eq:fproperties},
\begin{equation}\label{eq:cramergrowth}
  r_{k+1}\gg\frac{r_k^3}{(\log r_k)^4},\qquad k\geq1.
\end{equation}
All implied constants here and below may depend on the constant in \eqref{eq:cramer}, but not on the choice of $p$ once $p$ is sufficiently large.

Put $u_k=\log r_k$. By \eqref{eq:cramergrowth}, there exists $C_0>0$ such that
\begin{equation}\label{eq:logrec}
  u_{k+1}\geq 3u_k-4\log u_k-C_0.
\end{equation}
To retain the exact exponent at the endpoint, put $g(u)=3u-4\log u-C_0$ and $h(u)=u-2\log u-D$, where $D>0$ is chosen so that $2D\geq C_0+2\log 3$. For all sufficiently large $u$, $2u\leq g(u)\leq3u$, the function $h$ is increasing, and
\[
  h(g(u))-3h(u)
  =2D-C_0-2\log\frac{g(u)}{u}\geq0.
\]
By \eqref{eq:logrec}, $h(u_{k+1})\geq3h(u_k)$. Since $u_1\geq2\log(p-1)$ and $h(u_1)\geq u_1/2$ for sufficiently large $p$, iteration gives
\begin{equation}\label{eq:endpointgrowth}
  u_k\geq h(u_k)\geq 3^{k-1}h(u_1)\gg 3^k\log p.
\end{equation}

Consider the generalized prime system
\[
  \Pc_p=\big(\Pp\setminus\{p,r_1,r_2,\ldots\}\big)
           \cup\{\sqrt{r_1},\sqrt{r_2},\ldots\}.
\]
As in the proof of Proposition~\ref{prop:count}, distinct formal products have distinct numerical values. The prime $p$ is absent from $\Pc_p$. Replacing $B$ in \eqref{eq:exactcount} by
\[
  B_p(y)=\left\lfloor y\right\rfloor-\left\lfloor\frac{y}{p}\right\rfloor
        =\left(1-\frac1p\right)y+O(1),\qquad y\geq0,
\]
we obtain, for $\sqrt{r_k}\leq x<\sqrt{r_{k+1}}$,
\[
  N_{\Pc_p}(x)=\left(1-\frac1p\right)x
       \prod_{j=1}^k(1+r_j^{-1/2})+O(2^k).
\]
The constant in the estimate for $B_p$ is independent of $p$. By \eqref{eq:cramergrowth}, $r_{j+1}\geq4r_j$ for all $j\geq1$ when $p$ is sufficiently large. Applying \eqref{eq:pnormalization} and the tail estimate from the proof of Proposition~\ref{prop:count}, we obtain
\[
  |N_{\Pc_p}(x)-x|\ll 2^k+x\sum_{j>k}r_j^{-1/2}\ll2^k.
\]
By \eqref{eq:endpointgrowth} and $\log x\geq u_k/2$,
\[
  2^k=(3^k)^{\beta_*}\ll
       \frac{(\log x)^{\beta_*}}{(\log p)^{\beta_*}}.
\]
Thus, for some constant $C_1>0$ independent of $p$,
\[
  \limsup_{x\to\infty}\frac{|N_{\Pc_p}(x)-x|}{(\log x)^{\beta_*}}
  \leq \frac{C_1}{(\log p)^{\beta_*}}.
\]
Choose an ordinary prime $p$ sufficiently large that the right-hand side is strictly smaller than $c$. This proves \eqref{eq:smallconstant}.
\end{proof}

Theorem~\ref{thm:main} disproves Beurling's conjecture, whereas Proposition~\ref{prop:endpoint} identifies a smaller logarithmic order at which a conditional analogue of Olofsson's theorem holds. These results motivate the following conjecture.

\begin{conjecture}[Revised Beurling conjecture]
Let $\Pc$ be a generalized prime system. If
\begin{equation}\label{eq:revised}
  N_{\Pc}(x)=x+o\big((\log x)^{\log 2/\log 3}\big),\qquad x\to\infty,
\end{equation}
then $\Pc=\Pp$.
\end{conjecture}

Under Cram\'er's conjecture, Proposition~\ref{prop:endpoint} shows that the exponent in \eqref{eq:revised} cannot be replaced by any larger exponent, and that the little-$o$ condition cannot be replaced by a big-$O$ condition, even with an arbitrarily small positive constant. The proposition does not decide whether \eqref{eq:revised} itself implies $\Pc=\Pp$. In particular, it remains open whether the order of the error term in our construction can be further improved.

\end{document}